\documentclass[reqno]{amsart}
\usepackage{amssymb, amsmath, amsthm,latexsym,booktabs,todonotes, color, comment, lineno}

\usepackage[T1]{fontenc}
\usepackage[utf8]{inputenc}
\usepackage{lmodern}
\usepackage{amsmath,amssymb,amsthm}
\usepackage{geometry}
\usepackage{microtype}
\usepackage{listings}
\usepackage{xcolor}
\usepackage{hyperref}

\hypersetup{hidelinks}

\newcommand{\mex}{{\textnormal{mex}}}

\theoremstyle{definition}
\newtheorem{definition}{Definition}

\newtheorem{remark}[definition]{Remark}

\theoremstyle{plain}
\newtheorem{lemma}[definition]{Lemma}
\newtheorem{proposition}[definition]{Proposition}
\newtheorem{theorem}[definition]{Theorem}
\newtheorem{corollary}[definition]{Corollary}

\usepackage{todonotes}

\newcommand{\Nzero}{\mathbb{N}_0}
\newcommand{\Epass}{E_{\mathrm{pass}}}

\lstdefinestyle{pythonstyle}{
    language=Python,
    basicstyle=\ttfamily\footnotesize,
    keywordstyle=\color{blue!60!black},
    commentstyle=\color{gray!70!black},
    stringstyle=\color{green!40!black},
    showstringspaces=false,
    breaklines=true,
    keepspaces=true,
    columns=fullflexible,
    frame=single,
    framesep=5pt,
    xleftmargin=0pt,
    xrightmargin=0pt,
    tabsize=4
}

\begin{document}


\title[Expansions and rigidity in subtraction games with a single pass]{Compatible expansions and rigidity in\\ subtraction games with a single pass}
\author{Richard M. Low}
\address{Department of Mathematics\\
         San Jose State University\\
         San Jose, CA 95192\\
         USA}
\email{richard.low@sjsu.edu}


\date{November 29, 2026. Version 0.2}

\begin{abstract}
We study normal-play subtraction games with one globally
shared pass. An expansion of the subtraction set is
pass-equivalent if it leaves the pass-available
Sprague--Grundy sequence unchanged. We show that pure
periodicity of the ordinary sequence, together with an
eventual translation identity for the pass-available
sequence, forces every pass-equivalent expansion to
preserve the ordinary sequence. For two-element sets
$\{x,y\}$, we classify the compatible additional moves:
they are exactly $x+1,\ldots,y-1$ when $y\leq 2x$,
and there are none when $y>2x$. We extend this result
to sets in short intervals. An exact recurrence
certificate also exhibits a three-element set for
which pass-state rigidity fails.
\end{abstract}
\maketitle

\section{Introduction} \label{S:intro}

For a finite impartial game, the Sprague--Grundy value
of a position is the least nonnegative integer missing
from the values of its options. A position has value
zero precisely when it is losing for the next player
under optimal play. This numerical viewpoint
distinguishes game positions more finely than their
winning or losing status and makes it natural to ask
when a change of move rules preserves the entire
sequence of values.

Finite subtraction games are a basic setting for this
question. Their ordinary Grundy sequences are
eventually periodic, and the two-element case has an
explicit description \cite{Benesh}. The expansion
problem asks which new subtraction moves can be added
without changing the sequence; related work for
three-element sets appears in \cite{Ho}.

Let $S$ be a nonempty set of positive integers.
In an ordinary normal-play subtraction game, a move
removes $s\in S$ tokens from a heap of size $n$,
provided $s\leq n$. An expansion is a set $T\supseteq S$. Ordinary
expansion theory asks which additional moves preserve
the ordinary Sprague--Grundy sequence.

We consider a variant with one globally shared pass.
Either player may consume the pass without changing
the heap size, but the pass can be used only once
in total and is unavailable at terminal positions.

Write $g_S(n)$ for the ordinary Grundy value and
$h_S(n)$ for the value when the pass remains available.
The central rigidity question is whether
\[
h_T=h_S \quad\Longrightarrow\quad g_T=g_S
\]
for expansions $T\supseteq S$.

We establish three principal results.

First, a general sufficient condition for rigidity
follows from pure ordinary periodicity and an eventual
translation identity between $h_S$ and $g_S$.
Known two-move sequence formulas satisfy these
hypotheses, so every two-element subtraction set
is rigid.

Second, for $S=\{x,y\}$, $1\leq x<y$, we classify
all additional moves preserving the pass-state
sequence. They are precisely the integers strictly
between $x$ and $y$ when $y\leq 2x$, and no moves
otherwise. We also classify pass-equivalent expansions
of multi-move sets contained in sufficiently short
intervals.

Third, rigidity is not universal for three-element
bases. We establish the counterexample
\[
\{5,7,14\}\subset\{5,6,7,14\},
\]
whose pass-state sequences agree at every heap size
although the ordinary sequences differ. A finite
recurrence certificate proves the infinite sequence
identities.

The ordinary two-move Grundy words are known from
\cite{Benesh}, and the corresponding single-pass tail
identity is established in \cite{Low}. Related
pass-enabled subtraction games are considered in
\cite{Larsson}.

We first give the two-state recurrences and these
sequence formulas. We then establish a general
rigidity principle and work through the compatibility
and obstruction arguments before stating the complete
two-move classification. The final sections treat
short-interval sets and an explicit failure of
rigidity for a three-element base.

\section{Preliminaries}

For any finite set $B\subseteq\Nzero$, define
\[
\mex(B)=\min(\Nzero\setminus B).
\]

\noindent
For a nonempty subtraction set $S$, the ordinary
Grundy sequence satisfies
\begin{equation}
g_S(n)=
\mex\{g_S(n-s):s\in S,\ s\leq n\}.
\label{eq:ordinary}
\end{equation}

\noindent
Let $a=\min S$. Terminal positions satisfy
\[
g_S(n)=h_S(n)=0\qquad(0\leq n<a).
\]

\noindent
For $n\geq a$, the pass-state recurrence is
\begin{equation}
h_S(n)=
\mex\left(
\{h_S(n-s):s\in S,\ s\leq n\}
\cup\{g_S(n)\}
\right).
\label{eq:pass}
\end{equation}

\noindent
An expansion $T\supseteq S$ is
\emph{pass-equivalent} to $S$ if $h_T=h_S$.
We call $S$ \emph{pass-state rigid under expansion} if
\[
h_T=h_S\Longrightarrow g_T=g_S
\]
for every positive-move expansion $T\supseteq S$. We repeatedly use the elementary mex property
\begin{equation}
\mex(B\cup C)=\mex(B)
\quad\text{if}\quad \mex(B)\notin C.
\label{eq:mexproperty}
\end{equation}

\subsection*{The two-move sequences}

Fix
\[
S=\{x,y\},\qquad 1\leq x<y,
\]
and write
\[
g=g_S,\qquad h=h_S.
\]

\noindent
Put
\[
y=qx+r,\qquad q\geq1,\quad 0\leq r<x,
\qquad P=x+y.
\]

The following ordinary two-move period word is
known from \cite{Benesh}. Its explicit blocks are
the input for the compatibility arguments below.
The ordinary sequence is purely periodic.
Writing $a^j$ for $j$ copies of symbol $a$,
a period word is
\begin{equation}
G=
\begin{cases}
(0^x1^x)^{q/2}0^r2^{x-r}1^{r},
& q\text{ even},\\[2pt]
(0^x1^x)^{(q+1)/2}2^r,
& q\text{ odd},\ r>0,\\[2pt]
(0^x1^x)^{(q+1)/2},
& q\text{ odd},\ r=0.
\end{cases}
\label{eq:periodword}
\end{equation}

\noindent
The sequence $g$ is the infinite repetition of $G$.
The integer $P=x+y$ is always a period.

The pass-available sequence has a finite transient.
The classification in \cite{Low} supplies the eventual
translation identity
\begin{equation}
h(n)=g(n+y)\qquad(n\geq N),
\label{eq:tail}
\end{equation}
where
\begin{equation}
N=
\begin{cases}
(q+1)x,&q\text{ odd},\ r=0,\\
(q+2)x,&q\text{ even},\\
(q+2)x+r,&q\text{ odd},\ r>0.
\end{cases}
\label{eq:N}
\end{equation}

\noindent
In particular,
\begin{equation}
N\leq P+x.
\label{eq:Nbound}
\end{equation}

Define
\[
\Epass(x,y)=
\{z\in\mathbb Z_{>0}\setminus S:
h_{S\cup\{z\}}=h_S\}.
\]

\section{Technical lemmas and propositions}

We begin with a necessary condition for pass
equivalence. A translation of a purely periodic
ordinary word turns that condition into ordinary
sequence preservation. The two-move case then admits
a direct compatibility test.

The key observation is that pass-state equivalence
prevents an added move from connecting positions
with equal pass-state Grundy values. Under the
tail-translation hypothesis, pure periodicity
transfers this restriction to the ordinary sequence,
forcing ordinary equivalence.

\begin{lemma}
Suppose $T\supseteq S$ and $h_T=h_S$.
Every $z\in T\setminus S$ satisfies
\[
h_S(k+z)\neq h_S(k)\qquad(k\geq0).
\]
\end{lemma}

\begin{proof}
At heap size $k+z$, subtraction by $z$ is legal
and has pass-state option value $h_T(k)$.
The mex at $k+z$ cannot equal that option value.
Since $h_T=h_S$, the claimed inequality follows.
\end{proof}

\begin{theorem}
Let $S$ be a nonempty subtraction set.
Suppose $g_S$ is purely periodic with positive
period $P$. Assume that for some integers
$c,N\geq0$,
\begin{equation}
h_S(n)=g_S(n+c)\qquad(n\geq N).
\label{eq:generaltranslation}
\end{equation}
Then every expansion $T\supseteq S$ satisfies
\[
h_T=h_S\Longrightarrow g_T=g_S.
\]
The expansion $T$ may be infinite.
\end{theorem}

\begin{proof}
Suppose $h_T=h_S$. For every $z\in T\setminus S$,
Lemma~1 gives
\[
h_S(k+z)\neq h_S(k)\qquad(k\geq0).
\]
For $k\geq N$, the tail identity implies
\[
g_S(k+z+c)\neq g_S(k+c).
\]
Hence
\[
g_S(t+z)\neq g_S(t)\qquad(t\geq N+c).
\]

Fix $a\geq0$. Choose $m\geq0$ such that
$a+mP\geq N+c$. Pure periodicity gives
\[
g_S(a+z)=g_S(a+mP+z)
\neq g_S(a+mP)=g_S(a).
\]

Consequently, every additional legal ordinary
subtraction option avoids the original ordinary mex.
Strong induction using \eqref{eq:mexproperty} proves
$g_T=g_S$.

At any fixed heap size, only finitely many positive
moves are legal. Thus the argument also applies
when $T$ is infinite.
\end{proof}

\begin{corollary}
For $S=\{x,y\}$, $1\leq x<y$, and every expansion
$T\supseteq S$,
\[
h_T=h_S\Longrightarrow g_T=g_S.
\]
\end{corollary}

\begin{proof}
Apply Theorem~2 with $c=y$ and $P=x+y$.
\end{proof}

\begin{proposition}
For $z>0$, $z\notin S$,
\[
z\in\Epass(x,y)
\Longleftrightarrow
h(k+z)\neq h(k)\quad(k\geq0).
\]
Any collection of individually pass-compatible
moves may be adjoined simultaneously.
\end{proposition}

\begin{proof}
Necessity is Lemma~1.

Conversely, suppose $h(k+z)\neq h(k)$ for all
$k\geq0$. If $z<x$, then
$h(z)=h(0)=0$, a contradiction.

For $k\geq N$, the tail identity gives
\[
g(k+z+y)\neq g(k+y).
\]
Pure periodicity therefore yields
\[
g(k+z)\neq g(k)\qquad(k\geq0),
\]
as in Theorem~2.

Strong induction shows that every added ordinary
subtraction option avoids $g(n)$, so the ordinary
sequence remains unchanged.

The pass option consequently retains value $g(n)$.
Every additional pass-state subtraction option avoids
$h(n)$, so the pass-state mex also remains unchanged.

The argument applies simultaneously to any collection
of individually compatible moves.
\end{proof}

\begin{corollary}
Let $P=x+y$ and let $N$ be given by
\eqref{eq:N}. For $z>0$, $z\notin S$,
\[
z\in\Epass(x,y)
\Longleftrightarrow
h(k+z)\neq h(k)
\quad(0\leq k<N+P).
\]
\end{corollary}

\begin{proof}
The tail identity and ordinary periodicity give
\[
h(n+P)=h(n)\qquad(n\geq N).
\]
For $k\geq N$, both compared values belong to
the periodic tail. Every comparison therefore reduces
to one with $N\leq k<N+P$.
Proposition~4 completes the proof.
\end{proof}

\subsection{Compatible intermediate moves and exclusions}

We now apply the criterion to individual shifts.
The intermediate range is compatible in the
short-interval case; large shifts and short shifts
outside that range give the required obstructions.

Suppose
\[
x<y\leq2x,\qquad
y=x+r,\qquad 1\leq r\leq x.
\]
Then
\[
P=2x+r,\qquad g=0^x1^x2^r.
\]

\begin{lemma}
The residues $a\pmod P$ satisfying
\[
g(k+a)\neq g(k)\qquad(k\geq0)
\]
are precisely
\[
\{x,x+1,\ldots,x+r\}.
\]
\end{lemma}

\begin{proof}
If $0\leq a<x$, then $g(a)=g(0)=0$.

If $x+r<a<P$, take $k=P-a$.
Then $0<k<x$, and
\[
g(k+a)=g(P)=0=g(k).
\]

Conversely, let $a=x+u$, $0\leq u\leq r$.
In the cyclic word $0^x1^x2^r$, such a shift
sends each symbol to a different symbol.
Thus equal values never coincide.
\end{proof}

\begin{lemma}
For $y=x+r$, $1\leq r\leq x$,
\[
h[0,3x+r-1]=0^x2^x1^x3^r.
\]
The eventual tail begins at $N=3x+r$ with
periodic word $0^x1^x2^r$.
\end{lemma}

\begin{proof}
The first $x$ positions are terminal.

For $x\leq n<2x$, the subtraction options have
value zero and $g(n)=1$, so $h(n)=2$. For $2x\leq n<3x$, subtraction by $x$ reaches pass-state value $2$. A zero comes from subtraction
by $y$ when $n<2x+r$, and otherwise from the
ordinary pass option. Hence $h(n)=1$. For $3x\leq n<3x+r$, the option values are
$\{0,1,2\}$, so $h(n)=3$. Finally,
\[
N+y=4x+2r=2P,
\]
and \eqref{eq:tail} gives the claimed periodic tail.
\end{proof}

\begin{proposition}
If $x<z<y\leq2x$, then $z\in\Epass(x,y)$.
\end{proposition}

\begin{proof}
Write $z=x+u$, where $1\leq u<r$.
Ordinary noncoincidence follows from Lemma~6.

We check
\[
h(k+z)\neq h(k)
\]
for $0\leq k<N=3x+r$.

If $0\leq k<x$, then $h(k)=0$, while
$h(k+z)$ is $2$ or $1$.

If $x\leq k<2x$, then $h(k)=2$, while
$h(k+z)$ is $1$ or $3$.

If $2x\leq k<3x$, then $h(k)=1$.
When $k+z<N$, the shifted value is $3$.
Otherwise,
\[
0\leq k+z-N\leq x+u-r-1<x,
\]
so the shifted tail value is zero.

Finally, if $3x\leq k<N$, then $h(k)=3$,
while
\[
0\leq k+z-N\leq x+u-1<2x.
\]
Hence the shifted tail value is zero or one.

Thus pass-state noncoincidence holds for
$0\leq k<N$.

For $k\geq N$, the tail identity and ordinary
noncoincidence give
\[
h(k+z)=g(k+z+y)\neq g(k+y)=h(k).
\]
Proposition~4 completes the proof.
\end{proof}

\subsection{Large-move obstructions}

We first exclude moves near the ordinary period.

\begin{lemma}
Let $P=x+y$. No move
\[
y<z<P+x
\]
is pass-compatible.
\end{lemma}

\begin{proof}
If $y<z<P$, put $k=P-z$.
Then $0<k<x$, and
\[
g(k+z)=g(P)=0=g(k).
\]

If $P\leq z<P+x$, periodicity gives
\[
g(z)=g(z-P)=0=g(0).
\]

In either case ordinary noncoincidence fails.
By Corollary~3, the move cannot be pass-compatible.
\end{proof}

\begin{proposition}
If $y$ is not an even multiple of $x$,
no $z\geq N$ is pass-compatible.
\end{proposition}

\begin{proof}
Suppose $z\geq N$ is compatible.

For $0\leq k<x$,
\[
g(k)=h(k)=0.
\]
Pass-state compatibility gives
\[
h(z+k)\neq0.
\]

By Corollary~3, the ordinary sequence is also
preserved. Since subtraction by $z$ is a legal
ordinary option at $z+k$, the ordinary mex gives
\[
g(z+k)\neq g(k)=0.
\]

Moreover, the tail identity and $P=x+y$ yield
\[
h(z+k)=g(z+k+y)=g(z+k-x).
\]

Thus $g$ is nonzero throughout the $2x$
consecutive positions
\[
z-x\leq n<z+x.
\]

But the ordinary word \eqref{eq:periodword} has
maximum cyclic nonzero-run length
\[
\begin{cases}
x,&q\text{ even},\ r>0,\\
x+r,&q\text{ odd},\ r>0,\\
x,&q\text{ odd},\ r=0.
\end{cases}
\]

Every displayed length is less than $2x$.
These cases exhaust all parameters except even
multiples of $x$. The required nonzero run is
therefore impossible.

Thus a sufficiently large compatible move would
force an uninterrupted nonzero run of length $2x$
in the ordinary periodic word, which is impossible
outside the exceptional even-multiple cases
treated separately.
\end{proof}

\subsection{Exceptional even multiples}

Suppose
\[
y=2mx,\qquad m\geq2,
\]
and put
\[
L=2m+1,\qquad P=Lx.
\]

Then
\[
G=(0^x1^x)^m2^x,
\]
whose zero-block indices modulo $L$ are
\[
Z=\{0,2,\ldots,2m-2\}.
\]

\begin{lemma}
For $y=2mx$, $m\geq2$, the only residues
modulo $P$ satisfying ordinary noncoincidence
are $x$ and $y$.
\end{lemma}

\begin{proof}
The integer differences of elements of $Z$ form
\[
D=\{2j:-(m-1)\leq j\leq m-1\}.
\]

Since $L=2m+1$, these are the $L-2$ even
integers in
\[
[-(L-3),L-3].
\]

They are distinct modulo $L$: the difference
of two elements is even, has absolute value
less than $2L$, and cannot be a nonzero
multiple of the odd integer $L$.

Neither residue $1$ nor $-1$ occurs modulo $L$.
Indeed, representatives within the displayed
interval would have to be
$1$, $1-L$, $-1$, or $L-1$;
the odd representatives are not in $D$,
and the others lie outside the interval.

Consequently,
\[
Z-Z=(\mathbb Z/L\mathbb Z)
\setminus\{1,L-1\}.
\]

Write an arbitrary shift as
\[
z\equiv bx+u\pmod P,\qquad 0\leq u<x.
\]

If $u=0$, avoiding zero-to-zero collisions
requires $b\notin Z-Z$. Hence
\[
b\equiv1\text{ or }L-1\pmod L.
\]

If $0<u<x$, a shifted zero block overlaps
two consecutive blocks. Avoiding zero-to-zero
collisions would require both
\[
b\notin Z-Z,\qquad b+1\notin Z-Z.
\]
This is impossible because $1$ and $L-1$
are not consecutive modulo $L\geq5$.

Thus only the residues
\[
x,\qquad (L-1)x=y
\]
remain. Both are ordinarily noncoincident because
they correspond to original moves.
\end{proof}

\begin{proposition}
If $y$ is an even multiple of $x$, no additional
move $z\geq P=x+y$ is pass-compatible.
\end{proposition}

\begin{proof}
Write $y=2mx$, $m\geq1$. In every case,
\[
N=P+x.
\]

Lemma~9 excludes $P\leq z<N$.
Suppose $z\geq N$.

First assume $m\geq2$. Lemma~11 restricts
the residue of $z$ modulo $P$ to $x$ or $y$.

If $z\equiv x\pmod P$, then
\[
h(z)=g(z+y)=g(P)=0=h(0).
\]

If $z\equiv y\pmod P$, then
\[
h(z+x)=g(z+x+y)=g(y)=2=h(x).
\]

Both contradict necessary pass-state
noncoincidence.

Now let $m=1$, so $y=2x$ and $P=3x$.
Lemma~6 restricts an ordinarily noncoincident
residue $a=z\bmod P$ to
\[
x\leq a\leq2x.
\]

If $x\leq a<2x$, then
\[
h(z)=g(z+2x)=g(a+2x)=0=h(0).
\]

If $a=2x$, then
\[
h(z+x)=g(z+3x)=g(2x)=2=h(x).
\]

Again pass-state noncoincidence fails.
The proposition follows.
\end{proof}

\begin{theorem}
For every $1\leq x<y$,
\[
z\geq x+y
\Longrightarrow z\notin\Epass(x,y).
\]
\end{theorem}

\begin{proof}
Put $P=x+y$. Lemma~9 excludes
$P\leq z<P+x$.

For $z\geq P+x$, inequality
\eqref{eq:Nbound} gives $z\geq N$.
Proposition~10 handles all non-even-multiple
cases, and Proposition~12 handles every
even multiple.
\end{proof}

\subsection{Short moves}

\begin{lemma}
If $y>2x$, no additional move $z<y$
is pass-compatible.
\end{lemma}

\begin{proof}
For $0<z<x$,
\[
h(z)=h(0)=0.
\]

For $x<z<2x$, put $k=2x-z$.
Then $0<k<x$ and $2x<y$.
Since only subtraction by $x$ is legal below $y$,
\[
g(k)=g(2x)=0.
\]
Thus ordinary noncoincidence fails.

It remains to consider $2x\leq z<y$.
Write
\[
z=jx+t,\qquad
j=\lfloor z/x\rfloor,\qquad 0\leq t<x.
\]

Below $y$, only subtraction by $x$ is legal, so
\[
g(z)=j\bmod2.
\]

The pass-state recurrence gives initial block
values
\[
0,2,1,0,1,0,\ldots,
\]
indexed by $j=0,1,2,3,4,5,\ldots$.
Equivalently,
\[
h(jx+t)=
\begin{cases}
0,&j=0,\\
2,&j=1,\\
1,&j\geq2\text{ even},\\
0,&j\geq3\text{ odd}.
\end{cases}
\]

If $j$ is even, then $g(z)=g(0)=0$.
If $j$ is odd, then $j\geq3$ and
$h(z)=h(0)=0$.

Both alternatives exclude $z$.
\end{proof}

\section{Complete two-move classification}

The preceding propositions settle all possible
additional moves.

\begin{theorem}
For every $1\leq x<y$,
\[
\Epass(x,y)=
\begin{cases}
\{x+1,\ldots,y-1\},&y\leq2x,\\
\varnothing,&y>2x.
\end{cases}
\]
\end{theorem}

\begin{proof}
Suppose first that $y\leq2x$.

Proposition~8 shows that every $x<z<y$
is pass-compatible.

Moves $z<x$ fail pass-state noncoincidence.
Lemma~6 excludes all other additional moves
below $P$, while Theorem~13 excludes
$z\geq P$.

Hence
\[
\Epass(x,y)=\{x+1,\ldots,y-1\}.
\]

Now suppose $y>2x$.
Lemma~14 excludes all additional moves below $y$.
Lemma~9 excludes $y<z<P$, and
Theorem~13 excludes every $z\geq P$.

Therefore
\[
\Epass(x,y)=\varnothing.
\]
This completes the classification.
\end{proof}

\begin{remark}
For every integer $t\geq1$, the moves
\[
x+tP,\qquad y+tP
\]
preserve the ordinary Grundy sequence of $\{x,y\}$
by $P$-periodicity and mex induction.

However, Theorem~13 excludes all such moves
from the pass-compatible family.

Thus ordinary-compatible additions are infinite
in number, whereas pass-compatible additions
are finite, with cardinality
\[
\#\Epass(x,y)=
\begin{cases}
y-x-1,&y\leq2x,\\
0,&y>2x.
\end{cases}
\]
\end{remark}

\section{Short-interval expansions}

The classification determines all pass-equivalent
supersets of a two-move game.

\begin{theorem}
Let $A$ be a finite subtraction set with
\[
\min A=x,\qquad \max A=y,\qquad x<y\leq2x.
\]

Then
\[
g_A=g_{\{x,y\}},\qquad h_A=h_{\{x,y\}}.
\]

Moreover, for every expansion $T\supseteq A$,
\[
h_T=h_A
\Longleftrightarrow
A\subseteq T\subseteq\{x,x+1,\ldots,y\}.
\]

In particular, $A$ is pass-state rigid
under expansion.

The number of pass-equivalent expansions of $A$,
including $A$, is
\[
2^{\,y-x+1-|A|}.
\]

They form a Boolean lattice of rank
$y-x+1-|A|$.
\end{theorem}

\begin{proof}
The inclusions
\[
\{x,y\}\subseteq A\subseteq\{x,x+1,\ldots,y\}
\]
and Theorem~15 show that every additional move
of $A$ is individually pass-compatible
with $\{x,y\}$.

By Proposition~4, these moves may be adjoined
simultaneously. Thus
\[
h_A=h_{\{x,y\}},\qquad g_A=g_{\{x,y\}}.
\]

If $h_T=h_A$, then $h_T=h_{\{x,y\}}$.
Every additional move of $T$ must therefore
belong to the compatible interval, proving
\[
T\subseteq\{x,\ldots,y\}.
\]

Conversely, any such $T$ preserves the pass-state
sequence by simultaneous compatibility.

Two-move rigidity yields $g_T=g_A$ whenever
$h_T=h_A$.

Finally, the independently selectable additional
moves are precisely $\{x,\ldots,y\}\setminus A$,
whose cardinality is $y-x+1-|A|$.
This gives the count and Boolean-lattice structure.
\end{proof}

\section{A Three-move counterexample}

The two-move rigidity theorem does not extend
to every three-element subtraction set.

\begin{theorem}
Let
\[
S=\{5,7,14\},\qquad T=\{5,6,7,14\}.
\]
Then
\[
h_T=h_S,\qquad g_T\neq g_S.
\]
\end{theorem}

To establish equality for all heap sizes,
rather than for a finite sampled prefix,
we use an exact recurrence certificate.

\begin{lemma}
Let $S$ be a finite subtraction set with
maximum move $M$. Suppose candidate sequences
$\widetilde g,\widetilde h$ are periodic from
indices $0$ and $N$, respectively, with a
common positive period $P$.

If both candidates satisfy the ordinary and
pass-state mex recurrences at every index
\[
0\leq n<N+P+M,
\]
then
\[
\widetilde g=g_S,\qquad
\widetilde h=h_S
\]
at every nonnegative index.
\end{lemma}

\begin{proof}
The mex recurrences determine both sequences
uniquely. Thus the finite checks establish
the claimed sequences through the tested range.

For $n\geq N+P+M$ and every legal move $s\leq M$,
\[
n-s\geq N+P,\qquad n-P-s\geq N.
\]

\noindent
Both predecessor positions therefore lie in
the proposed periodic regime of $\widetilde h$.
The ordinary candidate $\widetilde g$ is
periodic from index zero. Hence all option values in the mex recurrences
at $n$ equal the corresponding option values
at $n-P$. The recurrence at $n-P$ has already been
checked or follows by induction. Therefore
the recurrence also holds at $n$. Induction establishes the candidate recurrences
at every subsequent index, and uniqueness
gives the result.
\end{proof}

\begin{proof}[Proof of Theorem~18]
The ordinary sequences have period words
\[
g_S=\mathtt{000001111122002331132},
\]
\[
g_T=\mathtt{000001111122002331134}.
\]

\noindent
Both words have length $21$.

The pass-state sequences have a common description
\[
h_S=h_T=H_0H_1^\omega,
\]
where
\[
H_0=\mathtt{00000222221111144220031133},
\]
\[
H_1=\mathtt{220000011122220031111}.
\]

\noindent
The prefix $H_0$ has length $26$, and the
period word $H_1$ has length $21$. Appendix~A provides a self-contained exact
checker for the ordinary and pass-state mex
recurrences at every index $0\leq n\leq60$.

Here $N=26$, $P=21$, and the maximum
subtraction move is $M=14$. Since
\[
N+P+M=61,
\]
Lemma~19 shows that the finite recurrence
certificate establishes the displayed
infinite sequences.

The pass-state descriptions coincide, but
\[
g_S(20)=2,\qquad g_T(20)=4.
\]
Therefore
\[
h_T=h_S,\qquad g_T\neq g_S.
\]

The example shows that equality of pass-available
Grundy sequences can conceal genuine differences
in the underlying ordinary sequences once the
subtraction set contains three moves.
\end{proof}

\begin{remark}
For any subtraction set $S$ and positive integer
$d$, let
\[
dS=\{ds:s\in S\}.
\]

The mex recurrences imply
\[
g_{dS}(dn+r)=g_S(n),
\]
\[
h_{dS}(dn+r)=h_S(n),
\qquad 0\leq r<d.
\]

Indeed, subtraction by $ds$ from $dn+r$
corresponds exactly to subtraction by $s$
from $n$, and terminal positions correspond.

Consequently, for every integer $d\geq1$,
the expansion
\[
\{5d,7d,14d\}\subset
\{5d,6d,7d,14d\}
\]
preserves the complete pass-state sequence
but not the ordinary sequence.

Thus failure of rigidity occurs in infinitely
many three-move base games.
\end{remark}

\section{Concluding remarks}

The general tail-translation theorem gives a
sufficient condition for pass-state rigidity
under expansion. Every two-move subtraction
game satisfies this condition.

For $S=\{x,y\}$, $1\leq x<y$,
the complete classification is
\[
\Epass(x,y)=
\begin{cases}
\{x+1,\ldots,y-1\},&y\leq2x,\\
\varnothing,&y>2x.
\end{cases}
\]

\noindent
This classification also determines
pass-equivalent expansions of finite multi-move
sets contained in sufficiently short intervals.

By contrast, the explicit three-move
counterexample shows that preservation of the
pass-state sequence does not universally force
preservation of the ordinary sequence. Thus rigidity is universal for two-move bases
and extends to broad multi-move families,
but it fails for general three-element
subtraction sets.

A natural open problem is to characterize
the finite subtraction sets that are
pass-state rigid under expansion.

\section{Tool and computational resource statement}\label{Disclosure}

The author supports the Leiden Declaration on Artificial Intelligence and Mathematics. GPT-6 Sol was used for exploratory investigations on the topic of this paper and for the creation of the Python 3 program in Appendix A.

\appendix

\section{Exact counterexample certificate}

The following Python 3 program checks both
mex recurrences for the periodic words used
in Theorem~18. No external Python packages
are required.

\vspace{10pt}

\begin{lstlisting}[style=pythonstyle]
def mex(values):
    used = set(values)
    k = 0
    while k in used:
        k += 1
    return k


H_PREFIX = "00000222221111144220031133"
H_PERIOD = "220000011122220031111"

G_S = "000001111122002331132"
G_T = "000001111122002331134"


def check_certificate(moves, g_word):
    period = 21
    start = 26
    max_move = max(moves)

    assert len(g_word) == period
    assert len(H_PREFIX) == start
    assert len(H_PERIOD) == period
    assert max_move == 14

    def g(n):
        return int(g_word[n % period])

    def h(n):
        if n < start:
            return int(H_PREFIX[n])
        return int(H_PERIOD[(n - start) % period])

    for n in range(start + period + max_move):
        legal = [s for s in moves if s <= n]

        expected_g = mex(
            g(n - s) for s in legal
        )
        assert g(n) == expected_g, (
            "g", moves, n
        )

        if n < min(moves):
            expected_h = 0
        else:
            options = [
                h(n - s) for s in legal
            ]
            options.append(g(n))
            expected_h = mex(options)

        assert h(n) == expected_h, (
            "h", moves, n
        )


check_certificate([5, 7, 14], G_S)
check_certificate([5, 6, 7, 14], G_T)

assert int(G_S[20]) == 2
assert int(G_T[20]) == 4

print("Periodic recurrence certificates passed.")
\end{lstlisting}

\vspace{10pt}

The certificate checks all indices
$0\leq n\leq60$. By Lemma~19, these checks
establish the infinite sequence identities.

\end{document}